\documentclass[a4paper, 12pt]{article}
\usepackage{amssymb, amsmath, amsfonts}
\usepackage[english,russian]{babel}
\usepackage{amsthm,amscd,amsfonts,latexsym,color,epsfig}

\begin{document}

\numberwithin{equation}{section}

\begin{center}
{\bf TOWARDS TO SOLUTION OF THE 
\\ FRACTIONAL TAKAGI -- TAUPIN EQUATIONS.\\ THE GREEN FUNCTION METHOD}
\end{center}
\begin{center}
{\bf  Murat O. Mamchuev $^1$, Felix N. Chukhovskii $^2$}
\end{center}


\small

{\bf Abstract.} 
Developing the comprehensive theory of the X-ray diffraction by distorted crystals remains to be topical of the mathematical physics. 
Up to now, the X-ray diffraction theory grounded on the Takagi -- Taupin equations with the first-order partial derivatives over the two coordinates within the X-ray scattering plane. 
In the work, the theoretical approach based on the first-order fractional Takagi -- Taupin equations with the 'quasi-time variable' of the order 
$\alpha\in(0,1]$ 
along the crystal depth has been suggested and the corresponding X-ray Cauchy issue is formulated. 
Accordingly, using the Green function method in the scope of the Cauchy issue, the fractional Takagi -- Taupin equations in the integral form have been derived. In the case of the inhomogeneous incident X-ray beam, the solution of the Cauchy issue of the X-ray diffraction by perfect crystal has been obtained and compared with the corresponding one based on the solution of the conventional Takagi -- Taupin equations, $\alpha=1.$ 
In turn, notice that the value of order $\alpha$ may be adjusted from the experimental X-ray diffraction data.
 \medskip

{\bf MSC 2010\/}: Primary 35F35;
                  
                  Secondary 35F40, 35A08, 35C15, 35L40, 45F05, 35Q70, 35Q92;

 \smallskip

{\bf Key Words and Phrases}: system of fractional partial differential
equations; Gerasimov -- Caputo fractional differentiation operator, X-ray diffraction tomography, transmission electron diffraction tomography, fractional Takagi -- Taupin equations; 

\normalsize

\section{Introduction}\label{sec:1}
\setcounter{equation}{0}

$\quad$ The dynamical theory of the X-ray diffraction by crystals has been based on the conventional Takagi -- Taupin (TT) equations 
\cite{Takagi-1962}, 
\cite{Takagi-1969}.
In some special cases, namely, in the case of the bent crystals with the constant deformation gradient of the reflection planes, analytical solutions of the Cauchy issue have been obtained 
\cite{Chukhovskii-1977}, \cite{Chukhovskii-1978}.

In the main, a number of the important physical results have been obtained due to the refined computer-numerical solutions of the TT-equations 
\cite{Epelboin-1983}, \cite{Honkanen-2018}. Noteworthy is the fact that progress in finding the TT-equations solutions is important aiming to facilitate and justify investigations of crystal-lattice defects in single crystals by using the X-ray diffraction (XRD) tomography technique. In present, the XRD tomography technique has got a good start to recover the 3D crystal-lattice defects key function $f\bf{(r) = hR(r)}$ involved in the TT-equations, $\bf{h}$ is the diﬀraction vector, $\bf{R(r)}$ is the displacement vector, ${\bf r}=(x,t)$ is the radius-vector within a the crystal sample (see, e.g., \cite{Chukhovskii-2020-1}, and references therein).

In recent decades, significant progress has been made in the development of the theory of equations with fractional derivatives. In particular, a theory of initial and boundary value problems for systems of partial differential equations of fractional order \cite{mamchuev-2010-du}--\cite{mamchuev-2017-fcaa} is constructed.
We will use these results to explore our model.

A goal of the present study is developing the theory of the X-ray diffraction based on the fractional TT-equations over of the 'quasi-time-variable' $t$ of the order $\alpha\in(0,1]$ along to the crystal depth. 
By using the Green function approach for solving the Cauchy problem, the integral fractional TT-equations are derived. 
In the case of the X-ray diffraction  by perfect crystals, the solution of the Cauchy problem has been obtained and compared with the suitable one of the conventional TT-equations for order $\alpha=1.$ 

\section{Fundamentals of the conventional diffraction theory}\label{sec:2}
\setcounter{equation}{0}

$\quad$ For complicity, let us start from the milestones of the theory of the X-ray diﬀraction. The Cauchy problem is formulated in the terms of the TT-equations (cf. \cite{Authier-2001}, \cite{Bowen-1998}, \cite{Chukhovskii-2020-1} for details).
To be specific, the latter can be cast in the matrix 
\begin{equation}\label{eq21}
{\bf S}(s_0,s_h){\bf E}(s_0,s_h)=i{\bf M}{\bf E}(s_0,s_h),
\end{equation}
$${\bf S}=\left (
\begin{array}{cc}
\frac{\partial}{\partial s_0}	& 0 \\
0  &   \frac{\partial}{\partial s_h}
\end{array}
\right ), \quad
{\bf M}=\sigma\left (
\begin{array}{cc}
0   & \exp[if({\bf r})]	 \\
\exp[-if({\bf r})]  &   0
\end{array}
\right ), $$
for the column vector   
${\bf E}(s_0,s_h)=\left (
\begin{array}{c}
E_0(s_0,s_h) \\
E_h(s_0,s_h)
\end{array}\right ),$ 
in search, where 
$E_0(s_0,s_h)$ and 
$E_h(s_0,s_h)$ are the refracted and diffracted amplitudes in the diffraction plane $(s_0,s_h),$ 
with the initial values on plane  $\tau=s_0+s_h,$  
$$
{\bf E}(x,t)\big|_{t=\tau}=\left (
\begin{array}{c}
E_0(x,\tau) \\
E_h(x,\tau)
\end{array}\right ).
$$

Hereafter is introduced the coupling coefficient $\sigma,$ 
$\sigma=-1+i\kappa,$ $0<\kappa<<1.$
In the matrix equation (\ref{eq21}), $f({\bf r})$ is the 3D phase function of no diagonal terms of matrix ${\bf M}$ and the dimensionless coordinates $s_0, s_h$ are linked with rectangular coordinates $x, t$ as follows $x = s_h-s_0$ and $t = s_h+s_0.$
Ultimately, putting on $\tau=0,$ the Cauchy problem is formulated as to searching the vector ${\bf E}(x,t)$ that is satisfied to the matrix equation 
\begin{equation}\label{eq23}
{\bf S}(x,t){\bf E}(x,t)=i{\bf M}(x,t){\bf E}(x,t),
\end{equation}
$$S=
\left(
\begin{array}{cc}
O_- & 0\\
0		 & O_+ 
\end{array}
\right), 
\quad O_{\pm}=\left(\frac{\partial}{\partial t}\pm\frac{\partial}{\partial x}\right);$$
$${\bf E}(x,0)=\left (
\begin{array}{c}
E_0(x,0) \\
E_h(x,0)
\end{array}\right ).$$
One step further, applying the double Fourier--Laplace $(k,p)$-transform to the matrix equation (\ref{eq23}), the Cauchy problem reduces to finding out vector ${\bf E}(x,t),$ which in turn satisfies to the integral matrix equation
$$\left(
\begin{array}{c}
E_0(k,p)\\
E_h(k,p)
\end{array}
\right)=\frac{1}{p^2+k^2+\sigma^2}
\left\{
\left(
\begin{array}{c}
(p+ik)E_0(k,0)\\
(p-ik)E_h(k,0)
\end{array}
\right)+
\right.$$
\begin{equation}\label{eq24}
\left.
+i\sigma
\left(
\begin{array}{l}
e^{if(x,0)}E_h(x,0)\\
e^{-if(x,0)}E_0(x,0)
\end{array}
\right)_k
+i\sigma
\left(
\begin{array}{l}
[O_+\, e^{if(x,t)}]\,E_h(x,t)\\
{[O_- e^{-if(x,t)}]} \,E_0(x,t)
\end{array}
\right)_{k,p}
\right\}.
\end{equation}
Hereafter the subscripts $k$ and $k,p$ denote the Fourier and Fourier--Laplace transforms, respectively.

Back using the Fourier--Laplace transform of (\ref{eq24}), one obtains 
\begin{equation}\label{eq25}
{\bf E}(x,t)=({\bf A}{\bf E}(x,t))(x,t)+({\bf B}{\bf E}(x,0))(x,t),
\end{equation}
where
$$({\bf A}{\bf E}(x,t))(x,t)=i\sigma
\int\limits_0^t dv
\int\limits_{x-t+v}^{x+t-v}
G(x-u,t-v)\times $$
$$\times
\left(
\begin{array}{cc}
0									&	{O_+\, e^{if(u,v)}} \\
{O_-\, e^{-if(u,v)}}	&		0
\end{array}
\right){\bf E}(u,v)du,$$
$$({\bf B}{\bf E}(x,0))(x,t)=
\int\limits_0^t dv
\int\limits_{x-t+v}^{x+t-v}
G(x-u,t-v)\times $$
$$\times
\left(
\begin{array}{cc}
O_+					&	i\sigma e^{if(u,0)} \\
i\sigma e^{-if(u,0)}	&		O_-
\end{array}
\right)\delta(v){\bf E}(u,0)du,$$
where 
$\delta(v)$ is the Dirac delta function and
the Green function  $G(x,t)$ is defined via following expression
$$G(x,t)=\frac{1}{i(2\pi)^2}\int\limits_{-i\infty}^{i\infty}dp
\int\limits_{-\infty}^{\infty}
\frac{e^{pt+ikx}}{p^2+k^2+\sigma^2}dk=$$
\begin{equation}\label{eq26}
=\frac{1}{2}J_0\left(\sigma \sqrt{t^2-x^2}\right)\Theta[t-|x|],
\end{equation}
here $J_0(z)$  is the zero-order Bessel
 function of argument   
$z.$  

Notice that in the basic integral equation (\ref{eq25}) integration is carried out over variables   within the triangle-shape area 
$$\Omega=\{(u,v): x-t+v<u<x+t-v, 0<v<t\}.$$  
Along with the conventional TT equations (\ref{eq21}), formulae (\ref{eq25}), (\ref{eq26}) represent by themselves integral formulation of the dynamical diffraction theory and they can be used for finding numerical solutions of the Cauchy diffraction problem  \cite{Epelboin-1983}, \cite{Honkanen-2018}.   At the same time, the integral formulation (\ref{eq25}), (\ref{eq26}) is preferred since it opens new ways for solving the Cauchy diffraction problem.

\section{The fractional Takagi--Taupin equations}\label{sec:2} 
\setcounter{equation}{0}

$\quad$ Let us introduce an additional parameter $\alpha\in(0,1]$ -- the order of the fractional time derivative, which can be verified based on experimental data.
Thus, as the basis of the model, we will consider a system of equations with partial derivatives of a fractional order not exceeding one, and we will call it the fractional Takagi -- Taupin equations (FTTEs).

One will consider the FTTEs in the matrix form
$$\left (
\begin{array}{cc}
\partial_{0t}^{\alpha}-\frac{\partial}{\partial x}	& 0  	 \\
 0   	&      \partial_{0t}^{\alpha}+\frac{\partial}{\partial x}
\end{array}
\right )  
\left (
\begin{array}{c}
E_0(x,t)	 \\
E_h(x,t)
\end{array}
\right ) =$$
\begin{equation}\label{eq31}
=i\sigma
\left (
\begin{array}{cc}
0						& \exp[if(x,t)]  	 \\
\exp[-if(x,t)]    	& 0      
\end{array}
\right )  
\left (
\begin{array}{c}
E_0(x,t)	 \\
E_h(x,t)
\end{array}
\right ),
\end{equation}
with the initial condition
\begin{equation}\label{eq32}
\left (
\begin{array}{c}
E_0(x,0)	 \\
E_h(x,0)
\end{array}
\right ) =
\left (
\begin{array}{c}
1	 \\
0
\end{array}
\right ), 
\quad -\infty<x<\infty.
\end{equation}
where 
$$
\partial_{0t}^{\alpha}g(t)=D_{0t}^{\alpha-1}\frac{d}{dt}g(t),
\quad 0<\alpha\leq 1, 
$$
is the Gerasimov -- Caputo fractional differentiation operator of order 
$\alpha$ \cite [p. 11]{n},
$D_{ay}^{\nu}$ is the Riemann -- Liouville fractional integro-differentiation operator
of order $\nu$  \cite [p. 9]{n}: 
$$D_{ay}^{\nu}g(y)=\frac{\mathop{\rm sgn}(y-a)}{\Gamma(-\nu)}
\int\limits_a^y \frac{g(s)ds}{|y-s|^{\nu+1}}, \quad  \nu< 0,$$
for $\nu \geq 0$ the operator $D_{ay}^{\nu}$ can be determined by 
recursive relation
$$D_{ay}^{\nu}g(y)=\mathop{\rm sgn}(y-a)\frac{d}{dy}D_{ay}^{\nu-1}g(y), \quad
\nu \geq 0,
$$
$\Gamma(z)$ is the Euler gamma-function.

Note that in the limit case $\alpha=1$ the operator $\partial_{0t}^{\alpha}g(t)$ go to the conventional partial derivatives 
$\frac{\partial}{\partial t}g(t).$


Acting on both sides of equation (\ref{eq31}) by the operator 
$diag(O_{+}^{\alpha}, O_{-}^{\alpha}),$ we obtain 
$$\left (
\begin{array}{cc}
O_{+}^{\alpha}	& 0  	 \\
 0   	&      O_{-}^{\alpha}
\end{array}
\right )  
\left (
\begin{array}{cc}
O_{-}^{\alpha}	& 0  	 \\
 0   	&     O_{+}^{\alpha}
\end{array}
\right )  
\left (
\begin{array}{c}
E_0(x,t)	 \\
E_h(x,t)
\end{array}
\right ) =$$
\begin{equation}\label{eq33}
=i\sigma
\left (
\begin{array}{cc}
0						& O_{+}^{\alpha}\exp[if(x,t)]E_h(x,t)  	 \\
O_{-}^{\alpha} \exp[-if(x,t)]E_0(x,t)    	& 0  
\end{array}
\right ),  
\end{equation}
where
$$O_{+}^{\alpha}=\partial_{0t}^{\alpha}+\frac{\partial}{\partial x},
\quad O_{-}^{\alpha}=\partial_{0t}^{\alpha}-\frac{\partial}{\partial x}.$$

After same transformation we get 
$$\left (
\begin{array}{cc}
O_{+}^{\alpha}O_{-}^{\alpha}	& 0  	 \\
 0   	&      O_{-}^{\alpha}O_{+}^{\alpha}
\end{array}
\right )  
\left (
\begin{array}{c}
E_0(x,t)	 \\
E_h(x,t)
\end{array}
\right ) =$$
$$
=
\left (
\begin{array}{cc}
0				& i\sigma e^{if}  	 \\
i\sigma e^{-if}     	& 0  
\end{array}
\right )
\left (
\begin{array}{cc}
O_{-}^{\alpha}	&     0 	 \\
0  	&     O_{+
}^{\alpha}
\end{array}
\right )  
\left (
\begin{array}{c}
E_0(x,t)	 \\
E_h(x,t)
\end{array}
\right ) 
 $$
$$+
\left (
\begin{array}{cc}
\partial_{0t}^{\alpha}				& 0  	 \\
0    	& \partial_{0t}^{\alpha}  
\end{array}
\right )
\left (
\begin{array}{cc}
0				& i\sigma e^{if}  	 \\
i\sigma e^{-if}     	& 0  
\end{array}
\right )
\left (
\begin{array}{c}
E_0(x,t)	 \\
E_h(x,t)
\end{array}
\right )+
$$
$$+
\left (
\begin{array}{cc}
0				& i\sigma \frac{\partial}{\partial x}e^{if}  	 \\
-i\sigma \frac{\partial}{\partial x}e^{-if}     	& 0  
\end{array}
\right )
\left (
\begin{array}{c}
E_0(x,t)	 \\
E_h(x,t)
\end{array}
\right )-$$
\begin{equation}\label{eq34}
-\left (
\begin{array}{cc}
0				& i\sigma e^{if}  	 \\
i\sigma e^{-if}     	& 0  
\end{array}
\right )
\left (
\begin{array}{cc}
\partial_{0t}^{\alpha}				& 0  	 \\
0    	& \partial_{0t}^{\alpha}  
\end{array}
\right )
\left (
\begin{array}{c}
E_0(x,t)	 \\
E_h(x,t)
\end{array}
\right ).  
\end{equation}

Using the fact that the vector $E(x,t)$ is a solution to equation (\ref{eq31}), we rewrite the equality (\ref{eq34}) in the form
$$\left (
\begin{array}{cc}
O_{+}^{\alpha}O_{-}^{\alpha}+\sigma^2	& 0  	 \\
 0   	&      O_{-}^{\alpha}O_{+}^{\alpha}+\sigma^2
\end{array}
\right )  
\left (
\begin{array}{c}
E_0(x,t)	 \\
E_h(x,t)
\end{array}
\right ) =$$
$$
=\left (
\begin{array}{cc}
\partial_{0t}^{\alpha}+if'_x				& 0  	 \\
0    	& \partial_{0t}^{\alpha}+if'_x  
\end{array}
\right )
\left (
\begin{array}{cc}
0				& i\sigma e^{if}  	 \\
i\sigma e^{-if}     	& 0  
\end{array}
\right )
\left (
\begin{array}{c}
E_0(x,t)	 \\
E_h(x,t)
\end{array}
\right )-
$$
\begin{equation}\label{eq35}
-\left (
\begin{array}{cc}
0				& i\sigma e^{if}  	 \\
i\sigma e^{-if}     	& 0  
\end{array}
\right )
\left (
\begin{array}{cc}
\partial_{0t}^{\alpha}				& 0  	 \\
0    	& \partial_{0t}^{\alpha}  
\end{array}
\right )
\left (
\begin{array}{c}
E_0(x,t)	 \\
E_h(x,t)
\end{array}
\right ).  
\end{equation}

Taking into account the next relation
$$[O_{\pm}^{\alpha}H(x,t)]_{k,p}=(p^{\alpha}\pm ik)[H(x,t)]_{k,p}-
p^{\alpha-1}[H(x,0)]_{k},$$ 
which follows from formula
\cite[p. 98]{Kilbas-2006}
$$[\partial_{0t}^{\alpha}H(x,t)]_{p}=p^{\alpha}[H(x,t)]_{p}-
p^{\alpha-1}H(x,0),$$ 
one can obtain
$$\left[O_{+}^{\alpha}O_{-}^{\alpha}E_0(x,t)\right]_{k,p}=
(p^{\alpha}+ ik)[O_{-}^{\alpha}E_0(x,t)]_{k,p}-
p^{\alpha-1}[O_{-}^{\alpha}E_0(x,0)]_{k}=$$
$$=(p^{\alpha}+ ik)(p^{\alpha}- ik)E_0(x,t)]_{k,p}-
p^{\alpha-1}(p^{\alpha}+ ik)[O_{-}^{\alpha}E_0(x,0)]_{k}-$$
\begin{equation} \label{eq36}
-p^{\alpha-1}[O_{-}^{\alpha}E_0(x,0)]_{k},
\end{equation}
$$\left[O_{-}^{\alpha}O_{+}^{\alpha}E_h(x,t)\right]_{k,p}=
(p^{\alpha}-ik)[O_{+}^{\alpha}E_h(x,t)]_{k,p}-
p^{\alpha-1}[O_{+}^{\alpha}E_h(x,0)]_{k}=$$
$$=(p^{\alpha}-ik)(p^{\alpha}+ik)E_h(x,t)]_{k,p}-
p^{\alpha-1}(p^{\alpha}-ik)[O_{+}^{\alpha}E_h(x,0)]_{k}-$$
\begin{equation} \label{eq37}
-p^{\alpha-1}[O_{+}^{\alpha}E_h(x,0)]_{k}.
\end{equation}
Further, taking into account the following relationships ({\it cf.} (\ref{eq31}))
$$O_{-}^{\alpha}E_0(x,0)=i\sigma e^{if(x,0)}E_h(x,0),$$
$$O_{+}^{\alpha}E_h(x,0)=i\sigma e^{-if(x,0)}E_0(x,0),$$
than from (\ref{eq35}) -- (\ref{eq37}) one obtains
$$\left (
\begin{array}{c}
E_0(x,t)	 \\
E_h(x,t)
\end{array}
\right )_{k,p}=
\frac{p^{\alpha-1}}{p^{2\alpha}+k^2+\sigma^2}
\left\{
\left (
\begin{array}{cc}
p^{\alpha}+ik				& 0  	 \\
0    	& p^{\alpha}-ik  
\end{array}
\right )
\left (
\begin{array}{c}
E_0(x,0)	 \\
E_h(x,0)
\end{array}
\right )_k+
\right.$$
$$
\left.
+i\sigma \left (
\begin{array}{c}
e^{if(x,0)}E_h(x,0)	 \\
e^{-if(x,0)}E_0(x,0)
\end{array}
\right )_k
\right\}+$$
$$+\frac{1}{p^{2\alpha}+k^2+\sigma^2}
\left\{
\left (
\begin{array}{cc}
\partial_{0t}^{\alpha}+if'_x				& 0  	 \\
0    	& \partial_{0t}^{\alpha}+if'_x  
\end{array}
\right )
\left (
\begin{array}{c}
i\sigma e^{if}E_h(x,t)	 \\
i\sigma e^{-if}E_0(x,t)
\end{array}
\right )-
\right. 
$$
\begin{equation}\label{eq38}
- 
\left. 
\left (
\begin{array}{cc}
0				& i\sigma e^{if}  	 \\
i\sigma e^{-if}     	& 0  
\end{array}
\right )
\left (
\begin{array}{c}
\partial_{0t}^{\alpha}E_0(x,t)	 \\
\partial_{0t}^{\alpha}E_h(x,t)
\end{array}
\right )
\right\}_{k,p}.  
\end{equation}
By applying the Efros theorem for operational calculus
\cite[p. 512]{Lavr},
the equality  \cite{Stankovic-1970}
$$\left(y^{\delta-1}\phi(-\beta,\delta;-ty^{-\beta})\right)_p=p^{-\delta}e^{-p^{\beta} t},$$
where
$$\phi(\alpha,\mu;z)=\sum\limits_{n=0}^{\infty}
\frac{z^n}{n!\Gamma(\alpha n+\mu)}$$
that is nothing else the Wright function \cite{Wright-1934},
%
%
one can see that the following expressions for the inversion of the Fourier--Laplace transform take place
$$\left(\frac{1}{p^{2\alpha}+k^2+\sigma^2}\right)_{x,t}=\frac{1}{i(2\pi)^2}\int\limits_{-i\infty}^{i\infty}dp\int\limits_{-\infty}^{\infty}
\frac{e^{pt+ikx}}{p^{2\alpha}+k^2+\sigma^2}dk=$$
$$=\frac{1}{2}\int\limits_{|x|}^{\infty}
J_0\left(\sigma\sqrt{\tau^2-x^2}\right) 
\frac{1}{t}\phi\left(-\alpha, 0;-\frac{\tau}{t^{\alpha}}\right)d\tau
=G_{\alpha}(x,t),$$

$$\left(\frac{p^{\alpha-1}}{p^{2\alpha}+k^2+\sigma^2}\right)_{x,t}=\frac{1}{i(2\pi)^2}\int\limits_{-i\infty}^{i\infty}dp
\int\limits_{-\infty}^{\infty}
\frac{p^{\alpha-1}e^{pt+ikx}}{p^{2\alpha}+k^2+\sigma^2}dk=
$$
$$=D_{0t}^{\alpha-1}G_{\alpha}(x,t)=
\frac{1}{2}\int\limits_{|x|}^{\infty}
J_0\left(\sigma\sqrt{\tau^2-x^2}\right) 
t^{-\alpha}\phi\left(-\alpha, 1-\alpha;-\frac{\tau}{t^{\alpha}}\right)d\tau.$$


Note that the  function $G_{\alpha}(x,t)$ properties have been investigated in the work \cite{mamchuev-2017-fcaa}.

Thus, one obtains the following integral matrix equation
\begin{equation}\label{eq39}
{\bf E}(x,t)=({\bf A}^{\alpha}{\bf E}(x,t))(x,t)+({\bf B}^{\alpha}{\bf E}(x,0))(x,t),
\end{equation}
where
$$({\bf A}^{\alpha}{\bf E}(x,t))(x,t)=-i\sigma\int\limits_{0}^{t}dv\int\limits_{-\infty}^{\infty}
G_{\alpha}(x-u,t-v)\cdot
\left\{D_{0v}^{\alpha-1}\frac{\partial}{\partial v}
[{\bf K}(u,v){\bf E}(u,v)]+\right.$$

\begin{equation}\label{eq310}
+\left.
\left (
\begin{array}{cc}
1		& 0  	 \\
0    	& -1  
\end{array}
\right )\left[\frac{\partial}{\partial u}{\bf K}(u,v)\right]
{\bf E}(u,v)-
{\bf K}(u,v)D_{0v}^{\alpha-1}\frac{\partial}{\partial v}
{\bf E}(u,v)\right\}du,
\end{equation}

$$({\bf B}^{\alpha}{\bf E}(x,0))=\int\limits_{0}^{t}dv\int\limits_{-\infty}^{\infty}
D_{tv}^{\alpha-1}G_{\alpha}(x-u,t-v)\times$$
\begin{equation}\label{eq311}
\times \left (
\begin{array}{cc}
O^{\alpha}_+				& i\sigma e^{if(u,0)}  	 \\
i\sigma e^{-if(u,0)}    	& O^{\alpha}_-  
\end{array}
\right )
{\bf E}(u,0)\delta(v)du,
\end{equation}

$${\bf K}(x,t)=
\left (
\begin{array}{cc}
0			  	&  e^{if(x,t)}  	 \\
e^{-if(x,t)}    & 0  
\end{array}
\right ). 
$$

Noteworthy is the fact that in the case when the key function $f(x,t)\equiv 0,$ the Cauchy problem in the integral matrix form (\ref{eq39})-(\ref{eq311}) has the rigorous solution as follows 
$$
{\bf E}(x,t)=({\bf B}^{\alpha}{\bf E}(x,0))(x,t),
$$
and/or in the explicit form 
$${\bf E}(x,t)=-\int\limits_{-\infty}^{\infty}D_{0t}^{2\alpha-1}
G_{\alpha}(x-u,t){\bf E}(u,0)du+$$
$$+\left (
\begin{array}{cc}
1		& 0  	 \\
0    	& -1  
\end{array}
\right )
\int\limits_{-\infty}^{\infty}D_{0t}^{\alpha-1}
G_{\alpha}(x-u,t){\bf E}'(u,0)du+$$
\begin{equation}\label{eq313}
+i\sigma
\left (
\begin{array}{cc}
0		& 1  	 \\
1    	& 0  
\end{array}
\right )
\int\limits_{-\infty}^{\infty}D_{0t}^{\alpha-1}
G_{\alpha}(x-u,t){\bf E}(u,0)du,
\end{equation}
respectively.

\section{The limited case of the fractional parameter $\alpha=1$}
\label{subsec:3.1}
\setcounter{equation}{0}

$\quad$ In the case when  $\alpha=1,$ the basic system (\ref{eq31}) reduces the form
(\ref{eq23}).

Let us pass in the relation (\ref{eq39}) to the limit at $\alpha\to 1.$

For this we need the following assertions.
 
{\bf Lemma 4.1.}
\cite{pizvran}. {\it Let the function $g(t)$ be absolutely integrable on any finite interval of the semiaxis $t>0,$ be continuous at the point $t=1$ and grow at $t\to \infty$ no faster than ${\rm exp}\{\sigma t^{\delta}\},$ $\sigma>0,$ $\delta<\frac{1}{1-\beta}.$
Then, }
$$\lim\limits_{\beta\to 1}\int\limits_{0}^{\infty}
g(t)\phi(-\beta,0;-t)dt=g(1), \quad
\lim\limits_{\beta\to 1}\int\limits_{0}^{\infty}
g(t)\phi(-\beta,\beta;-t)dt=\int\limits_{0}^{1}g(t)dt.$$

{\bf Lemma 4.2.}
{\it Under the conditions of Lemma 4.1, the relation
$$\lim\limits_{\beta\to 1}\int\limits_{0}^{\infty}
g(t)\phi(-\beta,1-\beta;-t)dt=g(1)$$
take place. }

Notice that Lemma 4.2 can be proved in the similar way as Lemma 4.1.

To be specific, the limit of the expression in the curly brackets in the right-hand side of (\ref{eq310}) is equal to
$$\left\{\frac{\partial}{\partial v}
[{\bf K}{\bf E}]+
\left (
\begin{array}{cc}
1		& 0  	 \\
0    	& -1  
\end{array}
\right )\left[\frac{\partial}{\partial u}{\bf K}\right]
{\bf E}-
{\bf K}\frac{\partial}{\partial v}
{\bf E}\right\}=
$$
$$=\left[\frac{\partial}{\partial v}{\bf K}
+\left (
\begin{array}{cc}
1		& 0  	 \\
0    	& -1  
\end{array}
\right )\frac{\partial}{\partial u}{\bf K}
\right]{\bf E}=
\left (
\begin{array}{cc}
O_+		& 0  	 \\
0    	& O_-  
\end{array}
\right ){\bf K}\cdot {\bf E}=$$
\begin{equation}\label{eq41}
=\left (
\begin{array}{cc}
0			  	&  [O_+e^{if(u,v)} ] 	 \\
{[O_-e^{-if(u,v)}]}    & 0  
\end{array}
\right ){\bf E}(u,v). 
\end{equation}

The next step, one goes to the limit as $\alpha\to 1$ for the
fundamental functions $G_{\alpha}(x,t)$ and $D_{0t}^{\alpha-1}G_{\alpha}(x,t).$
To do that, we rewrite them as follows
$$G_{\alpha}(x,t)=\frac{1}{2}\int\limits_{0}^{\infty}
J_0\left(\sigma\sqrt{(zt^{\alpha})^2-x^2}\right) 
t^{\alpha-1}\phi\left(-\alpha, 0;-z\right)\Theta(zt^{\alpha}-|x|)dz,$$
$$D_{0t}^{\alpha-1}G_{\alpha}(x,t)=\frac{1}{2}\int\limits_{0}^{\infty}
J_0\left(\sigma\sqrt{(zt^{\alpha})^2-x^2}\right) 
\phi\left(-\alpha, 1-\alpha;-z\right)\Theta(zt^{\alpha}-|x|)dz.$$

Going on to the limit $\alpha\to 1$ in the last expressions and taking into account Lemmas 4.1 and 4.2, one directly obtains
\begin{equation}\label{eq42}
\lim\limits_{\alpha \to 1}G_{\alpha}(x,t)=\frac{1}{2}
J_0\left(\sigma\sqrt{t^2-x^2}\right)\Theta(t-|x|)=G(x,t),
\end{equation}
\begin{equation}\label{eq43}
\lim\limits_{\alpha \to 1}D_{0t}^{\alpha-1}G_{\alpha}(x,t)=
\frac{1}{2}
J_0\left(\sigma\sqrt{t^2-x^2}\right)\Theta(t-|x|)=G(x,t).
\end{equation}

Thus, as if follows from (\ref{eq41})-(\ref{eq43}), in the limit case $\alpha=1,$ the basic integral equations (\ref{eq39}) reduces to (\ref{eq25}).

\section{\textsl{Quasi-Pendellösung} of the FTTEs}
\label{sec:1}
\setcounter{equation}{0}

As is known, in the case when $f(x,t)\equiv 0,$ the basic  TT-equations (\ref{eq23})  
with the initial conditions (\ref{eq32})
has the so-called Pendellösung (pendulum solutions) 
$$E_0=\cos (\sigma t), \quad E_h=i\sin (\sigma t).$$

Let us first consider a similar situation in the case of the FTTEs.
It is easy to show that in this case the system
$$
\left (
\begin{array}{cc}
\partial_{0t}^{\alpha}-\frac{\partial}{\partial x}	& 0  	 \\
 0   	&      \partial_{0t}^{\alpha}+\frac{\partial}{\partial x}
\end{array}
\right )  
\left (
\begin{array}{c}
E_0(x,t)	 \\
E_h(x,t)
\end{array}
\right ) = 
\left (
\begin{array}{cc}
0	& i\sigma  	 \\
i\sigma  & 0  
\end{array}
\right )  
\left (
\begin{array}{c}
E_0(x,t)	 \\
E_h(x,t)
\end{array}
\right ),
$$
with the initial condition
(\ref{eq32})
has the solution
\begin{equation} \label{eq52}
E_0(x,t)=E_{2\alpha, 1}(-\sigma^2t^{2\alpha}), \quad
E_h(x,t)=i\sigma t^{\alpha}E_{2\alpha, \alpha+1}(-\sigma^2t^{2\alpha}),
\end{equation}
where 
\begin{equation} \label{eq53}
E_{\rho, \mu}(z)=\sum\limits_{k=0}^{\infty}\frac{z^k}{\Gamma(\mu+\rho k)}, \quad
\rho>0, \, \mu>0
\end{equation}
is the  Mittag-Leffler type function
\cite[p. 117]{Djrb-1966}.

The validity of this statement follows from the next properties of the Mittag-Leffler type function
\begin{equation} \label{eq54}
E_{\alpha,\beta}(z)=\frac{1}{\Gamma(\beta)}+zE_{\alpha,\beta+\alpha}(z),
\end{equation}
\begin{equation} \label{eq55}
\partial_{0t}^{\mu}t^{\beta-1}E_{\alpha,\beta}(\lambda t^{\alpha})=
t^{\beta-\mu-1}E_{\alpha,\beta-\mu}(\lambda t^{\alpha}).
\end{equation}

Alternatively, it is clear, the same solution like (\ref{eq52}) can be directly obtained using the equation (\ref{eq313}) together with the initial conditions (\ref{eq32}).
Indeed,
$${\bf E}(x,t)=-
\left (
\begin{array}{c}
1	 	 \\
0    	  
\end{array}
\right )
\int\limits_{-\infty}^{\infty}D_{0t}^{2\alpha-1}
G_{\alpha}(x-u,t)du+$$
$$
+i\sigma
\left (
\begin{array}{cc}
0		& 1  	 \\
1    	& 0  
\end{array}
\right )
\left (
\begin{array}{c}
1	 	 \\
0    	  
\end{array}
\right )
\int\limits_{-\infty}^{\infty}D_{0t}^{\alpha-1}
G_{\alpha}(x-u,t)du=$$
\begin{equation}\label{eq56}
=\left (
\begin{array}{c}
D_{0t}^{\alpha}	 	 \\
i\sigma    	  
\end{array}
\right )
\int\limits_{-\infty}^{\infty}D_{0t}^{\alpha-1}
G_{\alpha}(x-u,t)du.
\end{equation}


Accordingly, by using the table integral \cite[p. 177]{Prudnikov-1983}
\begin{equation} \label{eq57}
\int\limits_{0}^{\tau}J_0(\sigma\sqrt{\tau^2-\xi^2})
d\xi=\int\limits_{0}^{\tau}\frac{z}{\sqrt{\tau^2-z^2}}J_0(\sigma z)
dz=\frac{\sin \sigma \tau}{\sigma}.
\end{equation}
and applying the Stancovic' transform \cite[p. 84]{PsMon}
$$\int\limits_{0}^{\infty}
t^{-\alpha}\phi\left(-\alpha,1-\alpha;-\frac{\tau}{t^{\alpha}}\right)\frac{\sin \sigma \tau}{\sigma}d\tau=
t^{\alpha}E_{2\alpha,\alpha+1}(-\sigma^2 t^{2\alpha}),$$
the straightforward routine calculations of integral in the right-hand side of equation (\ref{eq56}) yield
$${\bf E}(x,t)=\left (
\begin{array}{c}
D_{0t}^{\alpha}	 	 \\
i\sigma    	  
\end{array}
\right )
t^{\alpha}E_{2\alpha,\alpha+1}(-\sigma^2 t^{2\alpha})=$$
$$=\left (
\begin{array}{c}
t^{2\alpha}E_{2\alpha,1}(-\sigma^2 t^{2\alpha})	 	 \\
i\sigma t^{\alpha}E_{2\alpha,\alpha+1}(-\sigma^2 t^{2\alpha})   	  
\end{array}
\right ).$$

Note that for $\alpha=1$ from (\ref{eq53}) we obtain
$$E_{2\alpha}\left(-b^2t^{2\alpha}; 1\right)=E_{2}\left(-b^2t^{2}; 1\right)=
\sum\limits_{k=0}^{\infty}\frac{(-1)^k(bt)^{2k}}{\Gamma(1+2k)}=\cos(bt),$$
$$bt^{\alpha}E_{2\alpha}\left(-b^2t^{2\alpha}; \alpha+1\right)=btE_{2}\left(-b^2t^{2}; 2\right)=
\sum\limits_{k=0}^{\infty}\frac{(-1)^k(bt)^{2k+1}}{\Gamma(2+2k)}=\sin(bt).$$

The last relations give certain grounds for fractional generalization of the Takagi--Topen equations.

\section{Concluding remarks}

$\quad$ A goal of our study is to establish the mathematical framework for processing the reference 2D imaging patterns data of the X-ray diffraction tomography and then to develop mathematical background for solving the inverse tomography problem based on the general concept of the fractional Takagi--Taupin equations.	        
For the mathematical framework in processing of the reference 2D imaging patterns data, the integral matrix equation for solving the general Cauchy diffraction problem has been derived. 					    
Accordingly, it is shown that in the limit case when the fractional parameter $\alpha$ is equal to unity, all the results of solving the fractional Takagi--Taupin equations go on to the corresponding ones of the conventional X-ray diffraction theory including the known solutions in some particular cases, e.g., the pendulum one.         
One of the advantages of the present mathematical approach is the capacity numerically to incorporate the integral matrix equation to process it when the global minimum of the tomography target function in a $\chi^2$-sense 
needs to be achieved in a proper manner. 						   		 
Concluding, by using the mathematical framework and considering the results presented in this paper, we can claim that the fractional Takagi--Taupin equations approach is a good tool for obtaining digital structural crystal information from the reference 2D diffraction patterns tomography data. This would be a good topic for future research.






 \it
 
  \noindent
$^1$ Institute of Applied Mathematics and Automation, \\
Kabardino-Balkarian Scientific Center RAS, \\
``Shortanov" Str., 89A \\
360000  Nal'chik, Russian Federation \\[4pt]
e-mail: mamchuev@rambler.ru 
\\[10pt]

 \noindent
$^2$ A.V. Shubnikov Institute of Crystallography, 
\\
FSRC ``Cristallography and Photonics''  RAS,\\
Leninsky prospect, 59, \\
119333 Moscow,  Russian Federation  \\[4pt]
 e-mail: f\_chukhov@yahoo.ca 

\end{document}